\documentclass[reqno]{amsart}
\usepackage{amssymb,hyperref}

\theoremstyle{plain}
\newtheorem{thm}{Theorem}
\newtheorem{lem}{Lemma}
\newtheorem{cor}{Corollary}

\theoremstyle{definition}
\newtheorem{ex}{Example}

\renewcommand{\Re}{\mathrm{Re}}

\title
[Starlikeness and convexity for analytic functions ...]
{Starlikeness and convexity \\
for analytic functions \\
concerned with Jack's lemma}

\author{Hitoshi Shiraishi}
\address{Hitoshi Shiraishi \newline
Department of Mathematics \newline
Kinki University \newline
Higashi-Osaka, Osaka 577-8502, Japan}
\email{0733310104x@kindai.ac.jp}

\author{Shigeyoshi Owa}
\address{Shigeyoshi Owa \newline
Department of Mathematics \newline
Kinki University \newline
Higashi-Osaka, Osaka 577-8502, Japan}
\email{owa@math.kindai.ac.jp}

\subjclass[2000]{30C45}
\keywords{Analytic, univalent, starlike of order $\alpha$, convex of order $\alpha$.}

\date{}

\begin{document}

\begin{abstract}
There are many results for sufficient conditions of functions $f(z)$ which are analytic in the open unit disc $\mathbb{U}$ to be starlike and convex in $\mathbb{U}$.
The object of the present paper is to derive some interesting sufficient conditions for $f(z)$ to be starlike of order $\alpha$ and convex of order $\alpha$ in $\mathbb{U}$ concerned with Jack's lemma.
Some examples for our results are also considered with the help of Mathematica 5.2.
\end{abstract}

\begin{flushleft}
This paper was published in the journal: \\
Int. J. Open Probl. Comput. Sci. Math. {\bf 2} (2009), No. 1, 37--47. \\
\url{http://www.emis.de/journals/IJOPCM/files/IJOPCM(vol.2.1.3.M.9).pdf}
\end{flushleft}
\hrule

\

\

\maketitle

\section{Introduction}

\

Let $\mathcal{A}$ denote the class of functions $f(z)$ that are analytic in the open unit disk $\mathbb{U}=\{z \in \mathbb{C}:|z|<1\}$,
so that $f(0)=f'(0)-1=0$.

We denote by $\mathcal{S}$ the subclass of $\mathcal{A}$ consisting of univalent functions $f(z)$ in $\mathbb{U}$.
Let $\mathcal{S^{*}}(\alpha)$ be the subclass of $\mathcal{A}$ consisting of all functions $f(z)$ which satisfy
$$
\Re\left(\frac{zf'(z)}{f(z)}\right)>\alpha \qquad(z \in \mathbb{U})
$$
for some $ 0\leqq\alpha<1$.
A function $f(z) \in \mathcal{S^{*}}(\alpha)$ is sais to be starlike of order $\alpha$ in $\mathbb{U}$.
We denote by $\mathcal{S^{*}}=\mathcal{S^{*}}(0)$.

Also, let $\mathcal{K}(\alpha)$ denote the subclass of $\mathcal{A}$ consisting of functions $f(z)$ which satisfy
$$
\Re\left(1+\frac{zf''(z)}{f'(z)}\right)>\alpha  \qquad(z \in \mathbb{U})
$$
for some $0\leqq\alpha<1$.
A function $f(z)$ in $\mathcal{K}(\alpha)$ is said to be convex of order $\alpha$ in $\mathbb{U}$.
We say that $\mathcal{K}=\mathcal{K}(0)$.
From the definitions for $\mathcal{S^{*}}(\alpha)$ and $\mathcal{K}(\alpha)$,
we know that $f(z) \in \mathcal{K}(\alpha)$ if and only if $zf'(z) \in \mathcal{S^{*}}(\alpha)$.

Let $f(z)$ and $g(z)$ be analytic in $\mathbb{U}$.
Then $f(z)$ is said to be subordinate to $g(z)$ if there exists an analytic function $w(z)$ in $\mathbb{U}$ satisfying $w(0)=0$, $|w(z)| < 1$ $(z\in\mathbb{U})$ and $f(z)=g(w(z))$.

We denote this subordination by
$$
f(z) \prec g(z) \qquad ( z\in\mathbb{U}).
$$

\

The basic tool in proving our results is the following lemma due to Jack \cite{m1ref1} (also, due to Miller and Mocanu \cite{m1ref2}).

\

\begin{lem} \label{jack} \quad
Let $w(z)$ be analytic in $\mathbb{U}$ with $w(0)=0$.
Then if $\left|w(z)\right|$ attains its maximum value on the circle $\left|z\right|=r$ at a point $z_{0}\in\mathbb{U}$,
then we have $z_{0}w'(z_{0})=kw(z_{0})$, where $k\geqq1$ is a real number.
\end{lem}

\

\section{Main results}

\

Applying Lemma \ref{jack},
we drive the following result for the class $\mathcal{C}$.

\

\begin{thm} \label{m2thm1} \quad
If $f(z)\in\mathcal{A}$ satisfies
$$
\Re\left(1+\frac{zf''(z)}{f'(z)}\right)
< \frac{\beta+1}{2(\beta-1)}
\qquad( z\in\mathbb{U})
$$
for some real $2\leqq\beta<3$, or
$$
\Re\left(1+\frac{zf''(z)}{f'(z)}\right)
< \frac{5\beta-1}{2(\beta+1)}
\qquad (z\in\mathbb{U})
$$
for some real $1<\beta\leqq2$, then
$$
\frac{zf'(z)}{f(z)}
\prec \frac{\beta(1-z)}{\beta-z}
\qquad (z\in\mathbb{U})
$$
and
$$
\left|\frac{zf'(z)}{f(z)}-\frac{\beta}{\beta+1}\right|
< \frac{\beta}{\beta+1}
\qquad (z\in\mathbb{U}).
$$

This implies that $f(z) \in \mathcal{S}^*$ and $\displaystyle\int_{0}^{z}\dfrac{f(t)}{t}dt \in \mathcal{K}$.
\end{thm}

\

\begin{proof}\quad
Let us define the function $w(z)$ by
$$
\frac{zf'(z)}{f(z)}
= \frac{\beta(1-w(z))}{\beta-w(z)}
\qquad (w(z) \neq \beta).
$$

Clearly,
$w(z)$ is analytic in $\mathbb{U}$ and $w(0)=0$.
We want to prove that $|w(z)|<1$ in $\mathbb{U}$.
Since
$$
1+\frac{zf''(z)}{f'(z)}
= \frac{\beta(1-w(z))}{\beta-w(z)}-\frac{zw'(z)}{1-w(z)}+\frac{zw'(z)}{\beta-w(z)},
$$
we see that
\begin{align*}
\Re\left(1+\frac{zf''(z)}{f'(z)}\right)
&= \Re\left(\frac{\beta(1-w(z))}{\beta-w(z)}-\frac{zw'(z)}{1-w(z)}+\frac{zw'(z)}{\beta-w(z)}\right)\\
&< \frac{\beta+1}{2(\beta-1)}\qquad(z\in\mathbb{U})
\end{align*}
for $2\leqq\beta<3$ and
\begin{align*}
\Re\left(1+\frac{zf''(z)}{f'(z)}\right)
&= \Re\left(\frac{\beta(1-w(z))}{\beta-w(z)}-\frac{zw'(z)}{1-w(z)}+\frac{zw'(z)}{\beta-w(z)}\right)\\
&< \frac{5\beta-1}{2(\beta+1)}
\qquad (z\in\mathbb{U})
\end{align*}
for $1<\beta\leqq2$. If there exists a point $z_{0} \in \mathbb{U}$ such that
$$
\max_{\left| z \right| \leqq \left| z_{0} \right|} \left| w(z) \right|
= \left| w(z_{0}) \right|
= 1,
$$
then Lemma \ref{jack} gives us that $w(z_{0})=e^{i \theta}$ and $z_{0}w'(z_{0})=kw(z_{0})$ $(k \geqq 1)$.
Thus we have
\begin{align*}
1+\frac{z_{0}f''(z_{0})}{f'(z_{0})}
&= \frac{\beta(1-w(z_{0}))}{\beta-w(z_{0})}-\frac{z_{0}w'(z_{0})}{1-w(z_{0})}+\frac{z_{0}w'(z_{0})}{\beta-w(z_{0})}\\
&= \beta+\beta(1-\beta+k)\frac{1}{\beta-e^{i\theta}}-\frac{k}{1-e^{i\theta}}.
\end{align*}

If follows that
\begin{align*}
\Re\left(\frac{1}{\beta-w(z_{0})}\right)
&= \Re\left(\frac{1}{\beta-e^{i \theta}}\right)\\
&= \frac{1}{2\beta}+\frac{\beta^{2}-1}{2\beta(1+\beta^{2}-2\cos\theta)}
\end{align*}
and
\begin{align*}
\Re\left(\frac{1}{1-w(z_{0})}\right)
&= \Re\left(\frac{1}{1-e^{i \theta}}\right)\\
&= \frac{1}{2}.
\end{align*}

Therefore, we have
$$
\Re\left(1+\frac{z_{0}f''(z_{0})}{f'(z_{0})}\right)
= \frac{1+\beta}{2}+\frac{(\beta^{2}-1)(1-\beta+k)}{2(1+\beta^{2}-2\beta\cos\theta)}.
$$

This implies that, for $2 \leqq \beta<3$,
\begin{align*}
\Re\left(1+\frac{z_{0}f''(z_{0})}{f'(z_{0})}\right)
&\geqq \frac{1+\beta}{2}+\frac{(\beta+1)(1-\beta+k)}{2(\beta-1)}\\
&\geqq \frac{1+\beta}{2}+\frac{(\beta+1)(2-\beta)}{2(\beta-1)}\\
&= \frac{\beta+1}{2(\beta-1)}
\end{align*}
and, for $1 < \beta \leqq2$,
\begin{align*}
\Re\left(1+\frac{z_{0}f''(z_{0})}{f'(z_{0})}\right)
&\geqq \frac{1+\beta}{2}+\frac{(\beta-1)(1-\beta+k)}{2(\beta+1)}\\
&\geqq \frac{1+\beta}{2}+\frac{(\beta-1)(2-\beta)}{2(\beta+1)}\\
&= \frac{5\beta-1}{2(\beta+1)}.
\end{align*}

This contradicts the condition in the theorem.
Therefore, there is no $z_{0}\in\mathbb{U}$ such that $|w(z_{0})|=1$.
This means that $|w(z)| < 1$ for all $z\in\mathbb{U}$,
so that
$$
\frac{zf'(z)}{f(z)}
\prec \frac{\beta(1-z)}{\beta-z}\qquad(z\in\mathbb{U}).
$$

Furthermore,
since
$$
w(z)
= \frac{\beta\left(\dfrac{zf'(z)}{f(z)}-1\right)}{\dfrac{zf'(z)}{f(z)}-\beta}\qquad (z\in\mathbb{U})
$$
and $|w(z)|<1\ (z\in\mathbb{U})$,
we conclude that
$$
\left|\frac{zf'(z)}{f(z)}-\frac{\beta}{\beta+1}\right|
< \frac{\beta}{\beta+1}
\qquad (z\in\mathbb{U}),
$$
which implies that $f(z) \in \mathcal{S}^*$.
Furthermore,
we see that $f(z) \in \mathcal{S}^*$ if and only if $\displaystyle\int_{0}^{z}\dfrac{f(t)}{t}dt \in \mathcal{K}$.
\end{proof}

\

Taking $\beta=2$ in the theorem,
we have following corollary due to Singh and Singh \cite{m1ref3}.

\

\begin{cor} \label{m2cor1} \quad
If $f(z)\in\mathcal{A}$ satisfies
$$
\Re\left(1+\frac{zf''(z)}{f'(z)}\right)
< \frac{3}{2}\qquad (z\in\mathbb{U}),
$$
then
$$
\frac{zf'(z)}{f(z)}
\prec \frac{2(1-z)}{2-z}
\qquad (z\in\mathbb{U})
$$
and
$$
\left|\frac{zf'(z)}{f(z)}-\frac{2}{3}\right|
< \frac{3}{2}
\qquad (z\in\mathbb{U}).
$$
\end{cor}

\

With Theorem \ref{m2thm1},
we give the following example.

\

\begin{ex} \label{m2ex1} \quad
For $2 \leqq \beta < 3$,
we consider the function $f(z)$ given by
$$
f(z)=\frac{\beta-1}{2}\left(1-(1-z)^{\frac{2}{\beta-1}}\right)
\qquad (z \in \mathbb{U}).
$$

If follows that
$$
\frac{zf'(z)}{f(z)}=\frac{2z(1-z)^{\frac{3-\beta}{\beta-1}}}{(\beta-1)\left(1-(1-z)^\frac{2}{\beta-1}\right)}
\qquad(z \in \mathbb{U})
$$
and
\begin{align*}
\Re\left(1+\frac{zf''(z)}{f'(z)}\right)
&= \Re\left(\frac{\beta-1-2z}{(\beta-1)(1-z)}\right)\\
&= \Re\left(\frac{2}{\beta-1}-\frac{3-\beta}{(\beta-1)(1-z)}\right)\\
&< \frac{\beta+1}{2(\beta-1)}
\qquad (z \in \mathbb{U}).
\end{align*}

Therefore, the function $f(z)$ satisfies the condition in Theorem \ref{m2thm1}.
If we define the function $w(z)$ by
$$
\frac{zf'(z)}{f(z)}=\frac{\beta(1-w(z))}{\beta-w(z)}\qquad (w(z) \neq \beta),
$$
then we see that $w(z)$ is analytic in $\mathbb{U}$,
$w(0)=0$ and $|w(z)| < 1$ $(z \in \mathbb{U})$ with Mathematica 5.2.
This implies that
$$
\frac{zf'(z)}{f(z)} \prec \frac{\beta(1-z)}{\beta-z}\qquad (z \in \mathbb{U}).
$$

For $1 < \beta \leqq 2$, we consider
$$
f(z)=\frac{\beta+1}{2(2\beta-1)}\left(1-(1-z)^{\frac{2(2\beta-1)}{\beta+1}}\right)
\qquad (z\in \mathbb{U}).
$$

Then we have that
$$
\frac{zf'(z)}{f(z)}=\frac{2(2\beta-1)z(1-z)^{\frac{3(\beta-1)}{\beta+1}}}{(\beta+1)\left(1-(1-z)^{\frac{2(2\beta-1)}{\beta+1}}\right)}
$$
and
$$
\Re\left(1+\frac{zf''(z)}{f'(z)}\right)
=\Re\left(\frac{\beta+1-2(2\beta-1)z}{(\beta+1)(1-z)}\right)
< \frac{5\beta-1}{2(\beta+1)}
\qquad ( z \in \mathbb{U}).
$$

Thus,
the function $f(z)$ satisfies the condition in Theorem \ref{m2thm1}.
Define the function $w(z)$ by
$$
\frac{zf'(z)}{f(z)}=\frac{\beta(1-w(z))}{\beta-w(z)}\qquad (w(z) \neq \beta).
$$

Then $w(z)$ is analytic in $\mathbb{U}$,
$w(0)=0$ and $|w(z)| < 1$ $(z \in \mathbb{U})$ with Mathematica 5.2.
Therefore,
we have that
$$
\frac{zf'(z)}{f(z)} \prec \frac{\beta(1-z)}{\beta-z}
\qquad (z \in \mathbb{U}).
$$

In particular, if we take $\beta=2$ in this example,
then $f(z)$ becomes
$$
f(z)=z-\frac{1}{2}z^{2}\in\mathcal{S^{*}}.
$$
\end{ex}

\

\begin{thm} \label{m2thm2} \quad
If $f(z)\in\mathcal{A}$ satisfies
$$
\Re\left(1+\frac{zf''(z)}{f'(z)}\right)
> -\frac{\beta+1}{2\beta(\beta-1)}
\qquad(z\in\mathbb{U})
$$
for some real $\beta\leqq-1$, or
$$
\Re\left(1+\frac{zf''(z)}{f'(z)}\right)
> \frac{3\beta+1}{2\beta(\beta+1)}
\qquad(z\in\mathbb{U})
$$
for some real $\beta>1$, then
$$
\frac{f(z)}{zf'(z)}
\prec \frac{\beta(1-z)}{\beta-z}
\qquad(z\in\mathbb{U})
$$
and
$$
f(z)\in\mathcal{S^{*}}\left(\frac{\beta+1}{2\beta}\right).
$$

This implies that $\displaystyle\int_{0}^{z}\dfrac{f(t)}{t}dt \in \mathcal{K}\left(\frac{\beta+1}{2\beta}\right)$.
\end{thm}

\

\begin{proof}\quad
Let us define the function $w(z)$ by
\begin{equation}
\frac{f(z)}{zf'(z)}=\frac{\beta(1-w(z))}{\beta-w(z)}
\qquad (w(z)\neq\beta). \label{m2thm2eq1}
\end{equation}

Then, we have that $w(z)$ is analytic in $\mathbb{U}$ and $w(0)=0$.
We want to prove that $|w(z)|<1$ in $\mathbb{U}$.
Differentiating (\ref{m2thm2eq1}) in both side logarithmically and simplifying,
we obtain
$$
1+\frac{zf''(z)}{f'(z)}
= \frac{\beta-w(z)}{\beta(1-w(z))}+\frac{zw'(z)}{1-w(z)}-\frac{zw'(z)}{\beta-w(z)},
$$
and, hence
\begin{align*}
\Re\left(1+\frac{zf''(z)}{f'(z)}\right)
&= \Re\left(\frac{\beta-w(z)}{\beta(1-w(z))}+\frac{zw'(z)}{1-w(z)}-\frac{zw'(z)}{\beta-w(z)}\right)\\
&> -\frac{\beta+1}{2\beta(\beta-1)}\qquad (z\in\mathbb{U})
\end{align*}
for $\beta\leqq-1$, or
\begin{align*}
\Re\left(1+\frac{zf''(z)}{f'(z)}\right)
&= \Re\left(\frac{\beta-w(z)}{\beta(1-w(z))}+\frac{zw'(z)}{1-w(z)}-\frac{zw'(z)}{\beta-w(z)}\right)\\
&> \frac{3\beta+1}{2\beta(\beta+1)}\qquad (z\in\mathbb{U})
\end{align*}
for $\beta>1$.
If there exists a point $z_{0} \in \mathbb{U}$ such that
$$
\max_{\left| z \right| \leqq \left| z_{0} \right|} \left| w(z) \right|
= \left| w(z_{0}) \right| = 1,
$$
then Lemma \ref{jack} gives us that $w(z_{0})=e^{i \theta}$ and $z_{0}w'(z_{0})=kw(z_{0})$ $(k \geqq 1)$.
Thus we have
\begin{align*}
1+\frac{z_{0}f''(z_{0})}{f'(z_{0})}
&= \frac{\beta-w(z_{0})}{\beta(1-w(z_{0}))}+\frac{z_{0}w'(z_{0})}{1-w(z_{0})}-\frac{z_{0}w'(z_{0})}{\beta-w(z_{0})}\\
&= \frac{1}{\beta}+\frac{\beta-1}{\beta(1-e^{i\theta})}+\frac{k}{1-e^{i\theta}}-\frac{k\beta}{\beta-e^{i\theta}}.
\end{align*}

Therefore,
we have
$$
\Re\left(1+\frac{z_{0}f''(z_{0})}{f'(z_{0})}\right)
= \frac{1}{2}+\frac{1}{2\beta}-\frac{k(\beta^{2}-1)}{2(1+\beta^{2}-2\beta\cos\theta)}.
$$

This implies that,
for $\beta\leqq-1$,
\begin{align*}
\Re\left(1+\frac{z_{0}f''(z_{0})}{f'(z_{0})}\right)
&\leqq \frac{1}{2}+\frac{1}{2\beta}-\frac{k(\beta+1)}{2(\beta-1)}\\
&\leqq \frac{1}{2}+\frac{1}{2\beta}-\frac{\beta+1}{2(\beta-1)}\\
&= -\frac{\beta+1}{2\beta(\beta-1)}.
\end{align*}
and, for $\beta>1$,
\begin{align*}
\Re\left(1+\frac{z_{0}f''(z_{0})}{f'(z_{0})}\right)
&\leqq \frac{1}{2}+\frac{1}{2\beta}-\frac{k(\beta-1)}{2(\beta+1)}\\
&\leqq \frac{1}{2}+\frac{1}{2\beta}-\frac{\beta-1}{2(\beta+1)}\\
&= \frac{3\beta+1}{2\beta(\beta+1)}.
\end{align*}

This contradicts the condition in the theorem.
Therefore, there is no $z_{0}\in\mathbb{U}$ such that $|w(z_{0})|=1$.
This means that $|w(z)|<1$ for all $z\in\mathbb{U}$, this is, that
$$
\frac{f(z)}{zf'(z)}
\prec \frac{\beta(1-z)}{\beta-z}\qquad (z\in\mathbb{U}).
$$

Furthermore,
since
$$
w(z)
= \frac{\beta\left(1-\dfrac{zf'(z)}{f(z)}\right)}{1-\beta\dfrac{zf'(z)}{f(z)}}\qquad (z\in\mathbb{U})
$$
and $|w(z)|<1\ (z\in\mathbb{U})$,
we conclude that
$$
f(z)\in\mathcal{S^{*}}\left(\frac{\beta+1}{2\beta}\right).
$$

Noting that $f(z) \in \mathcal{S}^*(alpha)$ if and only if $\displaystyle\int_{0}^{z}\dfrac{f(t)}{t}dt \in \mathcal{K}(\alpha)$, we complete the proof of the theorem.
\end{proof}

\

For Theorem \ref{m2thm2},
we give the following example.

\

\begin{ex} \label{m2ex2} \quad
For $\beta>1$,
we take
$$
f(z)=\frac{\beta(\beta+1)}{-\beta^{2}+2\beta+1}\left(1-(1-z)^{\frac{-\beta^{2}+2\beta+1}{\beta(\beta+1)}}\right)
\qquad (z \in \mathbb{U}).
$$

Then, $f(z)$ satisfies
$$
\frac{zf'(z)}{f(z)}=\frac{(-\beta^{2}+2\beta+1)z}{\beta(\beta+1)(1-z)^{\frac{2\beta^{2}-\beta-1}{\beta(\beta+1)}}\left(1-(1-z)^{\frac{-\beta^{2}+2\beta+1}{\beta(\beta+1)}}\right)}
$$
and
\begin{align*}
\Re\left(1+\frac{zf''(z)}{f'(z)}\right)
&= \Re\left(\frac{\beta(\beta+1)+(\beta^{2}-2\beta-1)z}{\beta(\beta+1)(1-z)}\right)\\
&> \frac{3\beta+1}{2\beta(\beta+1)}
\qquad (z\in\mathbb{U}).
\end{align*}

Therefore,
$f(z)$ satisfies the condition of Theorem \ref{m2thm2}.
Let us define the function $w(z)$ by
$$
\frac{f(z)}{zf'(z)}=\frac{\beta(1-w(z))}{\beta-w(z)}
\qquad (w(z)\neq\beta).
$$

Then $w(z)$ is analytic in $\mathbb{U}$,
$w(0)=0$ and $|w(z)|<1\ z\in\mathbb{U}$ with Mathematica 5.2.
It follows that
$$
\frac{f(z)}{zf'(z)} \prec \frac{\beta(1-z)}{\beta-z}\qquad (z \in \mathbb{U}).
$$

Furthermore,
for $\beta\leqq-1$,
we consider the following function
$$
f(z)=-\frac{\beta(\beta-1)}{\beta^{2}+1}\left(1-(1-z)^{-\frac{\beta^{2}+1}{\beta(\beta-1)}}\right).
$$

Note that
$$
\frac{zf'(z)}{f(z)}=\frac{-(\beta^{2}+1)z}{\beta(\beta-1)(1-z)^{\frac{2\beta^{2}-\beta+1}{\beta(\beta-1)}}\left(1-(1-z)^{-\frac{\beta^{2}+1}{\beta(\beta-1)}}\right)}
$$
and
\begin{align*}
\Re\left(1+\frac{zf''(z)}{f'(z)}\right)
&= \Re\left(\frac{\beta(\beta-1)+(\beta^{2}+1)z}{\beta(\beta-1)(1-z)}\right)\\
&> \frac{\beta+1}{2\beta(\beta-1)}
\qquad (z\in\mathbb{U}).
\end{align*}

This implies that $f(z)$ satisfies the condition of Theorem \ref{m2thm2}.
Definning the function $w(z)$ by
$$
\frac{f(z)}{zf'(z)}=\frac{\beta(1-w(z))}{\beta-w(z)}\qquad(w(z)\neq\beta),
$$
we see that $w(z)$ is analytic in $\mathbb{U}$,
$w(0)=0$ and $|w(z)|<1\ (z\in\mathbb{U})$ with Mathematica 5.2.

Thus we have that
$$
\frac{f(z)}{zf'(z)} \prec \frac{\beta(1-z)}{\beta-z}
\qquad(z \in \mathbb{U}).
$$

Making $\beta=-1$ for $f(z)$,
we have
$$
f(z)=\frac{z}{1-z}\in\mathcal{K}.
$$
\end{ex}

\


\begin{thebibliography}{}

\bibitem{m1ref1}
I. S. Jack,
{\it Functions starlike and convex of order $\alpha$},
J. London Math. Soc. {\bf 3}(1971), 469--474.

\bibitem{m1ref2}
S. S. Miller and P. T. Mocanu,
{\it Second-order differential inequalities in the complex plane},
J. Math. Anal. Appl. {\bf 65}(1978), 289--305.

\bibitem{m1ref3}
R. Singh and S. Singh,
{\it Some sufficient conditions for univalence and starlikeness},
Coll. Math. {\bf 47}(1982), 309--314.

\end{thebibliography}
\end{document}